\def\draft{n}
\documentclass{amsart}
\usepackage[headings]{fullpage}
\usepackage{amssymb,epic,eepic,epsfig,amsbsy,amsmath}

\theoremstyle{plain}

\newtheorem{theorem}{Theorem}
\newtheorem{proposition}{Proposition}[section]

\newtheorem{corollary}[proposition]{Corollary}

\theoremstyle{definition}
\newtheorem{definition}[proposition]{Definition}

\theoremstyle{remark}

\newtheorem{remark}[proposition]{Remark}

\def\printname#1{
        \if\draft y
                \smash{\makebox[0pt]{\hspace{-0.5in}
                        \raisebox{8pt}{\tt\tiny #1}}}
        \fi
}

\newcommand{\psdraw}[2]
         {\begin{array}{c} \hspace{-1.3mm}
        \raisebox{-4pt}{\epsfig{figure=draws/#1.eps,width=#2}}
        \hspace{-1.9mm}\end{array}}

\newlength{\standardunitlength}
\catcode`\@=11
\long\def\@makecaption#1#2{%
     \vskip 10pt

\setbox\@tempboxa\hbox{
       \small\sf{\bfcaptionfont #1. }\ignorespaces #2}%
     \ifdim \wd\@tempboxa >\captionwidth {%
         \rightskip=\@captionmargin\leftskip=\@captionmargin
         \unhbox\@tempboxa\par}%
       \else
         \hbox to\hsize{\hfil\box\@tempboxa\hfil}%
     \fi}
\font\bfcaptionfont=cmssbx10 scaled \magstephalf
\newdimen\@captionmargin\@captionmargin=2\parindent
\newdimen\captionwidth\captionwidth=\hsize
\catcode`\@=12

\def\lbl#1{\label{#1}\printname{#1}}

\def\BN{\mathbb N}
\def\BZ{\mathbb Z}

\def\BQ{\mathbb Q}
\def\BR{\mathbb R}
\def\BC{\mathbb C}

\def\a{\alpha}

\def\Ga{\Gamma}

\def\ga{\gamma}

\def\e{\epsilon}
\def\Ga{\Gamma}

\def\th{\theta}

\def\Li{\mathrm{Li}}

\def\Res{\mathrm{Res}}
\def\Li{\mathrm{Li}}

\def\ep{\epsilon}

\renewcommand\Re{\mathrm{Re}}

\begin{document}


\title[Resurgence of the fractional polylogarithms]{
Resurgence of the fractional polylogarithms}

\author{Ovidiu Costin}
\address{Department of Mathematics \\
         Ohio State University \\ 
         231 W 18th Avenue \\
         Columbus, OH 43210, USA \\ \newline
         {\tt http://www.math.ohio-state.edu/$\sim$costin } }
\email{costin@math.ohio-state.edu}
\author{Stavros Garoufalidis}
\address{School of Mathematics \\
         Georgia Institute of Technology \\
         Atlanta, GA 30332-0160, USA \\ 
         {\tt http://www.math.gatech} \newline {\tt .edu/$\sim$stavros } }
\email{stavros@math.gatech.edu}

\thanks{O.C. was supported in part by NSF grants DMS-0406193 and 
DMS-0600369 and S.G was supported in part by NSF grant DMS-0505445. \\
\newline
1991 {\em Mathematics Classification.} Primary 57N10. Secondary 57M25.
}

\date{July 16, 2009} 


\begin{abstract}
The fractional polylogarithms, depending on a complex parameter $\a$,
are defined by a series which is analytic
inside the unit disk. After an elementary conversion of the series
into an integral presentation, we show that the fractional polylogarithms
are multivalued analytic functions in the complex plane minus $0$ and $1$. 
For non-integer values of $\a$, we prove the analytic continuation, 
compute the monodromy around $0$ and $1$, give
a Mittag-Leffler decomposition and compute the asymptotic behavior for
large values of the complex variable. The fractional polylogarithms are
building blocks of resurgent functions that are used in proving 
that certain power series associated with knotted objects are
resurgent. This is explained in a separate publication \cite{CG3}.
The motivic or physical interpretation of the monodromy of the
fractional polylogarithms for non-integer values of $\a$ is unknown
to the authors.
\end{abstract}

\maketitle

\tableofcontents

\section{Introduction}
\lbl{sec.intro}

\subsection{The fractional polylogarithm and its history}
\lbl{sub.formalpseries}

For a complex number $\a$, let us define the {\em $\a$-polylogarithm} 
function $\Li_{\a}(z)$ by the following series:

\begin{equation}
\lbl{eq.polylog}
\Li_\a(z)=
\sum_{n=1}^\infty \frac{z^n}{n^{\a}} 
\end{equation}
which is absolutely convergent for $|z| <1$. These functions appear in
algebraic geometry, number theory, mathematical physics, applied mathematics
and the theory of special functions. Since

\begin{equation}
\lbl{eq.derLi}
z \frac{d}{dz} \Li_{\a}(z)= \Li_{\a-1}(z)
\end{equation}
we really need to study $\Li_{\a}(z)$ for $\a \bmod \BZ$. 

For integer $\a$ a lot is known about the $\a$-polylogarithm. For example, 
$\Li_0(z)=1/(1-z)$, thus \eqref{eq.derLi} implies that
for all $\a \in \BZ^-$, $\Li_{\a}(z) \in \BQ(z)$ is a rational 
function with a single singularity at $z=1$.

When $\a \in \BN$, the functions $\Li_{\a}(z)$ were studied in the nineteenth
century, forgotten for many years, and  rediscovered by
the algebraic geometers in the late 1970s; see for example Lewin's
book \cite{Lw}, Bloch's paper \cite{Bl} and the survey
articles \cite{Oe,Za1,Za2}. It is well known that $\Li_{\a}(z)$ is a multivalued 
function defined on $\BC\setminus\{0,1\}$ with computable monodromy; 
see \cite{We,Ha,BD, Oe} and \cite{MPV}. For $\a \in \BN$, the 
$\a$-polylogarithms are special functions that play a key role in
algebraic geometry. For $\a \in \BN$, the special values 
\begin{equation}
\lbl{zetaL}
\Li_{\a}(1)=\zeta(\a)
\end{equation}
are well-known examples of {\em periods}; see \cite{KZ}. 
This is not an accident. Zagier and Deligne conjectured that 
special values (at integers)
of $L$-functions of algebraic varieties are expressed
by the $\a$-polylogarithm for $\a \in \BN$; see \cite{Za1} and \cite{De}.
A motivic interpretation
of $\Li_{\a}(z)$ for $\a \in \BN$ is given in \cite{BD}, as well as a 
conjecture that the $\a$-th Beilinson-Deligne regulator maps are 
expressed by the $\a$-polylogarithm for $\a \in \BN$. 

For integer $\a$, {\em elliptic polylogarithms} that resemble $\Li_{\a}(e^z)$
were introduced by Beilinson-Levin in \cite{BL}, and further studied in
\cite{Lv} in relation to motivic cohomology conjectures for elliptic
curves. For a recent survey on the better-known dilogarithm, see \cite{Za2}.

The $\a$-polylogarithms for noninteger $\a$ are also classical and modern
objects. They were studied in the eighteenth
century by Jonqui\`ere as a function of two complex variables $\a$ and $z$;
see \cite{Jo}. Several approximation formulas were obtained by Jonqui\`ere
and half a century later by Truesdell, whose motivation was asymptotic
properties of {\em polymer structures}; see \cite{Tr}. Further results 
regarding approximation and analytic continuation were obtained by Pickard in 
the sixties, and more recently by Kenyon-Wilson in relation to {\em resonance}
of some {\em statistical mechanical models}; see \cite{Pi,KW} and also 
\cite[Prop.1]{CLZ}.

The $\a$-polylogarithm functions for half-integer $\a$
appear naturally in the context of an {\em Euler-MacLaurin summation},
and are also used in proving resurgence of some power series associated
to knotted objects; see \cite{CG2} and \cite{CG3}. They also 
play a prominent role in proving analytic continuation of some power 
series that encode quantum invariants of knotted objects; see for example 
\cite[Sec.7]{Ga}.

In addition, in 1994, M. Kruskal proposed to the first author 
to study the analytic continuation and the global bahavior of the function 
$\Li_{1/2}(z)$. This problem was a 
motivation for a global reconstruction theorem of resurgent functions from 
local data, developed by the first author several years ago (and independently
by \'Ecalle in \cite{Ec2}), and recently written in \cite{C}.

The purpose of this short note is to study the 
\begin{itemize}
\item[(a)]
the analytic continuation
\item[(b)]
the Mittag-Leffler decomposition
\item[(c)]
the asymptotic behavior for large $|z|$
\end{itemize}
of the polylogarithm function $\Li_{\a}(z)$
for non-integer $\a$. With over a century of 
history on the fractional polylogarithm, some of our results resemble 
statements of the above mentioned literature. However, we were not able to 
find the key Equation \eqref{eq.ML}, nor an explicit computation of the 
monodromy around $z=0$ and $z=1$ in the literature. The latter does not 
seem to have a finite dimensional faithful representation, and its motivic 
or physical origin is unknown when $\a \in \BQ\setminus\BZ$.

\subsection{Statement of the results}
\lbl{sub.results}

Let us recall first what is a multivalued analytic function on
$\BC\setminus\{0,1\}$. Such functions are examples of 
{\em global analytic functions} (see \cite{Ah}) and examples
of {\em resurgent functions} in the sense of \'Ecalle, \cite{Ec1}.
Let $X$ denote the
universal cover of $(\BC\setminus\{0,1\},1/2)$ with base point at $1/2$.
As a set, we have:
\begin{equation}
X=\left\{ \text{homotopy classes $[c]$ of paths $c$ in} 
\,\, \BC\setminus\{0,1\} \,\,
\text{starting at} \,\, \frac{1}{2} \right\}.
\end{equation}
There is an action of $F=\pi_1(\BC\setminus\{0,1\},1/2)$ on $X$ given
by $g \cdot [c]=[\gamma.c]$ for $g=[\a\ga]$ and $[c] \in X$. By a 
{\em multivalued analytic} (in short, {\em resurgent}) function
$f$ on $\BC\setminus\{0,1\}$ we mean an analytic
function on $X$. For $[c] \in X$, where $c$ is a path from $1/2$ to $z$,
we write, following \cite{Oe}:
\begin{equation}
\lbl{eq.analcont}
f^{[c]}(z):=f([c]).
\end{equation}  
Observe that $F$ is a free group on $[c_0]$ and $[c_1]$, where 
$$
c_0(t)=\frac{1}{2} e^{2 \pi i t}, \qquad
c_1(t)=1+\frac{1}{2} e^{2 \pi i t}
$$
are paths around $0$ and $1$ respectively:
$$
\psdraw{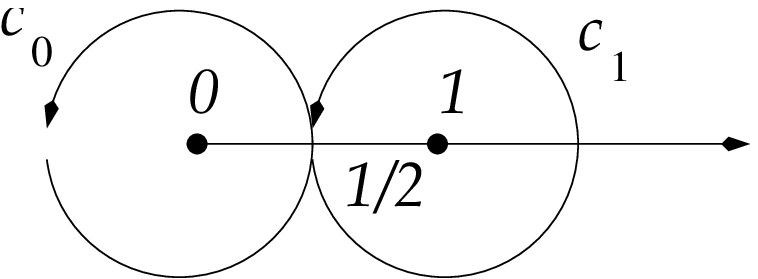}{1.5in}
$$

In what follows, $\a$ is {\em not} an integer. Let us introduce some
useful notation. 
Let $\ga$ denote a {\em Hankel contour} that encircles the positive real 
axis:
$$
\psdraw{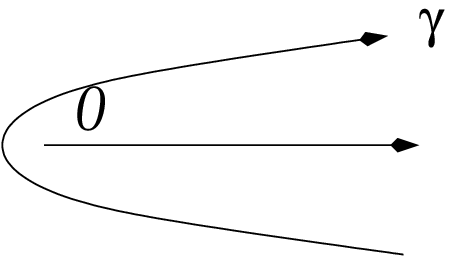}{1in} 
$$
The next definition uses notation familiar to algebraic geometry. See
for example the survey paper \cite{Oe}.

\begin{definition}
\lbl{def.M}
For $\a\in\BC\setminus\BZ$, let 
$M_{\a}(z)$ denote the multivalued function on $\BC\setminus\{0,1\}$
given by:
\begin{equation}
\lbl{eq.Ma}
M_{\a}(z)=
C_{\a} \, (\log z)^{\a-1} 
\end{equation}
where 
\begin{equation}
\lbl{eq.Ca}
C_{\a}=e^{\pi i(-\a-1)} \Ga(1-\a)
\end{equation}
For $k \in \BZ$, let us define the twisted multivalued functions 
$M_{\a}[k](z)$ for $z \in \BC\setminus\{0,1\}$ by:
\begin{equation}
\lbl{eq.twistedM}
M_{\a}[k](z):=M_{\a}(z \, e^{2 \pi i k})= 
C_{\a} \, (\log z + 2 \pi i k)^{\a-1}.
\end{equation}
\end{definition}

The following theorem converts the series 
\eqref{eq.polylog} of $\Li_{\a}(z)$ into an integral,
from which we can easily deduce the existence of analytic continuation.

\begin{theorem}
\lbl{thm.1}
\rm{(a)}
For $|z| <1$ and $\a$ such that $\Re(\a) >0$, 
$\Li_{\a}(z)$ has an  integral representation:
\begin{eqnarray}
\lbl{eq.Llaplace}
\Li_{\a}(z) &=& 
\frac{1}{\Ga(\a)} \int_0^\infty q^{\a-1} \frac{z}{e^q-z}
dq 
\end{eqnarray}
known as {\em Appell's integral} in \cite[Sec.2]{Tr}.
\newline
\rm{(b)} 
For $|z| <1$ and $\a \in \BC\setminus\BZ$, $\Li_{\a}(z)$ has an
 integral representation:
\begin{eqnarray}
\lbl{eq.Lhankel}
\Li_{\a}(z)  &=& \frac{C_{\a}}{2 \pi i} \int_{\ga} q^{\a-1}  \frac{z}{
e^q-z} dq
\end{eqnarray}
\rm{(c)} For all $\a \in \BC\setminus\BZ$,
$\Li_{\a}(z)$ has an analytic continuation to a multivalued function on 
$\BC\setminus\{0,1\}$. More precisely, let $z \in \BC\setminus\{0,1\}$
and $c$ any path from $1/2$ to $z$ in $\BC\setminus\{0,1\}$. Then,
we have:
\begin{align}
\lbl{eq.m0}
\Li_{\a}^{[c_0c]}(z) &= \Li_{\a}^{[c]}(z) &
\Li_{\a}^{[c_1c]}(z) &= \Li_{\a}^{[c]}(z)-(1-e^{2 \pi i \a})M_{\a}^{[c]}(z) \\
\lbl{eq.m1}
M_{\a}^{[c_0c]}(z) &= M_{\a}^{[c]}[1](z) &
M_{\a}^{[c_1c]}(z) &= -(1-e^{2 \pi i \a}) M_{\a}^{[c]}(z)
\end{align}
\rm{(d)} For $\a$ such that $\Re(\a) <0$,
$\Li_{\a}(z)$ has a Mittag-Leffler type decomposition:
\begin{equation}
\lbl{eq.ML}
\Li_{\a}(z) = 
C_{\a} \left( (\log z)^{\a-1}+\sum_{k=1}^\infty 
(\log z+ 2 \pi i k )^{\a-1} + (\log z- 2 \pi i k )^{\a-1} \right)
\end{equation}
where the series is uniformly convergent on compact sets. Thus, 
we have:
\begin{equation}
\lbl{eq.ML2}
\Li_{\a} = \sum_{k \in \BZ}  M_{\a}[k] :=\lim_{N \to \infty} \sum_{k=-N}^N  M_{\a}[k].
\end{equation}
\end{theorem}
When $\a$ is a negative integer, the right hand side of \eqref{eq.ML2}
is an {\em Eisenstein series}; see \cite{Ap}. 
The Mittag-Leffler decomposition 
\eqref{eq.ML2} is an analogue of {\em Hurwitz's theorem}; see \cite{Ap}.
The Mittag-Leffler \eqref{eq.ML2} implies is the following corollary.

\begin{corollary}
\lbl{cor.zetas}
For $\a$ such that $\Re(\a)<0$ and $z$ such that $\Re(z)<0$ and
$|z|<2\pi$ we have:
\begin{equation}
\lbl{eq.zetas}
\Li_{\a}(e^z)=C_{\a} z^{\a-1} + \sum_{n=0}^\infty \frac{\zeta(\a-n)}{n!} z^n.
\end{equation}
\end{corollary}
Compare with \cite[Prop.1]{CLZ}.

The integral formula \eqref{eq.Llaplace} and some stationary phase 
implies the following estimate for the behavior of the fractional
polylogarithms for large $|z|$.

\begin{corollary}
\lbl{cor.largez}
For $\Re(\a)>0$ and $z$ large we have:
\begin{equation}
\lbl{eq.largez}
\Li_{\a}(z)=-\frac{1}{\Gamma(\a+1)}\left((\log z)^{\a}+o(1)\right).
\end{equation}
\end{corollary}
For $\a \in \BN$, this is known; see \cite[Eqn.7]{Oe}.

\subsection{Plan of the proof}
\lbl{sub.plan}

Once we convert the series definition of the $\a$-polylogarithms into
an integral formula, analytic continuation follows from a general 
principle,
i.e., by moving the contour of integration and achieving analytic 
continuation. If we move the contour of integration to $-\infty$,
and the integral vanishes at $-\infty$, collecting residues gives
a Mittag-Leffler type decomposition of $\Li_{\a}(z)$ for $\a<0$, 
$\a\notin\BZ$.

\subsection{Acknowledgement}
An early version of this paper was presented at talks in Orsay
and the University of Maryland in the fall of 2006. 
The authors wish to thank J. \'Ecalle for encouraging conversations.
M. Kontsevich pointed out to the second author that some aspects
of the fractional polylogarithms have been studied independently by 
M. Kontsevich and D. Zagier. After the paper was written, J. Morava 
informed us of \cite{EM}, where the fractional polylogarithms are
also studied from the point of view of distributions over the real numbers.

\section{Proofs}
\lbl{sec.thm1}

\subsection{Proof of Theorem \ref{thm.1}}
\lbl{sub.thm1}

In this section we give a proof of Theorem \ref{thm.1}.

\begin{proof}
(a) For $\Re(\a)>0$ and $n \in \BN^+$ we have:
$$
\frac{1}{n^{\a}}=\frac{1}{\Ga(\a)} \int_0^\infty q^{\a-1}
e^{-n q} dq
$$
Interchanging summation and integration (valid for $|z|<1$) gives:

\begin{eqnarray*}
\Li_{\a}(z) &=& \sum_{n=1}^\infty \frac{1}{n^{\a}} z^n \\
&=&
\frac{1}{\Ga(\a)} \int_0^\infty q^{\a-1} \sum_{n=1}^\infty
(ze^{-q})^n dq \\
&=& 
\frac{1}{\Ga(\a)} \int_0^\infty q^{\a-1} \frac{z}{e^q-z} dq.
\end{eqnarray*}

\noindent
(b) Let 
\begin{equation}
\lbl{eq.Ia}
I_{\a}(z)=C_{\a} \, \int_{\ga} q^{\a-1}  \frac{z}{
e^q-z} dq
\end{equation}
denote the right hand side of Equation \eqref{eq.Lhankel}. Observe that
$I_{\a}(z)$ is well-defined for $\a \in \BC\setminus\BZ$, and 
$z\in\BC\setminus[1,\infty)$. 

Since for fixed $z$ inside the unit disk, 
both sides of \eqref{eq.Lhankel} are analytic functions of $\a \in
\BC\setminus\BZ$, it suffices to prove \eqref{eq.Lhankel} for 
$\Re(\a) >0$, $\a\not\in\BZ$.  
We claim that for such $\a$, we have:
\begin{equation}
\lbl{eq.link}
\int_0^\infty q^{\a-1} \frac{z}{e^q-z} dq=
\frac{1}{1-e^{2 \pi i \a}} 
\int_{\ga} q^{\a-1}  \frac{z}{e^q-z} dq.
\end{equation}
Indeed, we push the Hankel contour $\ga$ until its upper (resp. lower) 
part touches $\BR^+$ from above (resp. below) and push the tip of
the contour to touch zero.
On the upper part we have $q^{\a-1}=|q|^{\a-1}$, and on the lower side we have
\begin{equation}
\lbl{eq.jumpqa}
q^{\a-1}=e^{(\a-1)\ln|q|+2 \pi i (\a-1)}=|q|^{\a-1}e^{ 2\pi i \a}.
\end{equation} 
Moreover, the upper integral is traversed in the direction $(0,\infty)$ 
while the lower one is traversed from $(\infty,0)$. We thus get
\begin{equation}
\int_{\ga} q^{\a-1}\frac{1}{e^q-z} dq=
(1-e^{2\pi i \a})\int_0^{\infty}q^{\a-1}\frac{1}{e^q-z} dq
\end{equation}
Thus, \eqref{eq.link} follows. Since for $\a \in\BC\setminus\BZ$ the
$\Ga$ function satisfies the the reflexion symmetry (see eg. \cite{Ol}):

\begin{equation}
\lbl{eq.doubleGamma}
\frac{1}{(1-e^{2 \pi i \a}) \Ga(\a)}= \frac{e^{\pi i (-\a-1)} \Ga(1-\a)}{2 \pi i}
=\frac{C_{\a}}{2 \pi i} \, ,
\end{equation}
(b) follows. 

\noindent
(c) Fix $\a$ such that $\Re(\a)>0$, $\a \notin\BZ$. 
The integral representation \eqref{eq.Llaplace} analytically
continues $\Li_{\a}(z)$ for $z$ in the cut plane 
$\BC\setminus [1,\infty)$. Let us compute the {\em variation}
(i.e., jump) of the function across the cut $z \in (1,\infty)$. Changing
variable to $e^q=x$ in \eqref{eq.Llaplace}, we have:
$$
\Li_{\a}(z)=\frac{1}{\Ga(\a)} \int_1^\infty \frac{(\log x)^{\a-1}}{x}
\frac{z}{x-z} dx
$$
Fix $z \in (1,\infty)$. Then the above equality gives by contour
deformation and Cauchy's theorem (see eg. \cite{Df}):
\begin{eqnarray*}
\lim_{\ep \to 0^+}(\Li_{\a}(z+i \ep)-\Li_{\a}(z-i \ep)) &=&
\frac{2 \pi i}{\Ga(\a)} \, \Res\left( 
\frac{(\log x)^{\a-1}}{x}\frac{z}{x-z},x=z \right) \\
&=&  \frac{2 \pi i}{\Ga(\a)} (\log z)^{\a-1} \\
&=&(1-e^{2 \pi i\a}) M_{\a}(z). 
\end{eqnarray*}
On the other hand, Equation \eqref{eq.jumpqa} implies that
\begin{eqnarray*}
\lim_{\ep \to 0^+}(M_{\a}(z+i \ep)-M_{\a}(z-i \ep)) &=&
(1-e^{2 \pi i\a}) M_{\a}(z).
\end{eqnarray*}
Since
\begin{equation}
\lbl{eq.fc1}
f^{[c_1 c]}(z)=-\lim_{\e\to 0^+}( f^{[c]}(z+i \e)-f^{[c]}(z-i \e)),
\end{equation}
the above equations imply that
\begin{eqnarray*}
\Li_{\a}^{[c_1c]}(z) &=& \Li_{\a}^{[c]}(z)-(1-e^{2 \pi i \a})M_{\a}^{[c]}(z) \\
M_{\a}^{[c_1c]}(z) &=& -(1-e^{2 \pi i \a}) M_{\a}^{[c]}(z).
\end{eqnarray*}
On the other hand, \eqref{eq.Llaplace} defines an analytic function
for $z \in (-\infty,0)$, and the monodromy of $M_{\a}(z)$ for $z=0$
can be computed from the definition of $M_{\a}(z)$. Thus, we obtain
\begin{eqnarray*}
\Li_{\a}^{[c_0c]}(z) &=& \Li_{\a}^{[c]}(z) \\
M_{\a}^{[c_0c]}(z) &=& M_{\a}^{[c]}[1](z).
\end{eqnarray*}
This proves that when $\Re(\a)>0$, $\a\notin \BZ$,
$\Li_{\a}(z)$ is a multivalued function on 
$\BC\setminus\{0,1\}$ with monodromy given by \eqref{eq.m0} and 
\eqref{eq.m1}.
If $\a\in\BC\setminus\BZ$, use \eqref{eq.derLi}, the fact
\begin{equation}
\lbl{eq.derM}
z \frac{d}{dz} M_{\a}(z)= M_{\a-1}(z)
\end{equation}
and differentiation to conclude (c).

(d) Since both sides of \eqref{eq.ML} are analytic functions of $\a$
for fixed $z$, it suffices to prove \eqref{eq.ML} for 
$\a$ such that $\Re(\a)<0$, $\a\notin\BZ$. For such $\a$, we will use the 
integral representation of $\Li_{\a}(z)$ given by \eqref{eq.Lhankel}. 
Fix a complex number $z \in \BC\setminus[1,\infty)$
and the Hankel contour $\ga$ which separates the plane into two regions
so that $2 \pi i k + \log z$ lies in the region that contains $-\infty$
for all $k \in \BZ$. This is possible since the points $ 2 \pi i k + \log z$
lie in a vertical line. Now, push the Hankel contour to the left, and 
deform it to $-\infty$. Since $\Re(\a)<0$, the integral vanishes when 
the contour is deformed to $-\infty$. 
In the process of deformation, we apply Cauchy's
theorem and collect residues at the singularities.
The singularities of the integrand are simple poles at the points where
$q = \log z + 2 \pi i k$ for integer $k$. The residue is given by:
$$
\Res\left(q^{\a-1} \frac{z}{e^q-z}, q= \log z + 2 \pi i k\right)=
(\log z + 2 \pi i k)^{\a-1}
$$
When we push the contour to $-\infty$, we collect the series \eqref{eq.ML}
which is absolutely convergent on compact sets.
The result follows. 
\end{proof}

\begin{remark}
\lbl{rem.balt}
Part (b) of Theorem \ref{thm.1} states that for all $\a\in\BC\setminus\BZ$
and $|z|<1$ we have:
$$
\Li_{\a}(z)=I_{\a}(z).
$$
Moreover, $\Li_{\a}(z)$ satisfies the differential equation \eqref{eq.derLi}.
It is easy to show independently from Theorem \ref{thm.1} 
that for every $\a\in\BC\setminus\BZ$,
$I_{\a}(z)$ satisfies the differential equation
\begin{equation}
\lbl{eq.derI}
z \frac{d}{dz} I_{\a}(z)=I_{\a-1}(z).
\end{equation}
Indeed, use the algebraic identity:
\begin{equation}
\lbl{eq.partials}
\frac{d}{dz} \frac{z}{e^q-z}=\frac{e^q}{(e^q-z)^2}=
\frac{d}{dq} \frac{-1}{e^q-z}
\end{equation}
After differentiation and integration by parts, we have:
\begin{eqnarray*}
z I_{\a}'(z) &=& C_{\a}  z \int_{\ga}
q^{\a-1}  \frac{e^q}{(e^q-z)^2} dq \\
&=& C_{\a}  z \int_{\ga}
q^{\a-1} \frac{d}{dq} \frac{-1}{e^q-z}  dq \\
&=& C_{\a} (\a-1)  z \int_{\ga}
q^{\a-2}  \frac{1}{e^q-z}  dq \\
&=& C_{\a-1}  z \int_{\ga}
q^{\a-2}  \frac{1}{e^q-z}  dq \\
&=& I_{\a-1}(z).
\end{eqnarray*}
\end{remark}

\begin{remark}
\lbl{rem.MLexp}
An alternative way to prove part (d) of Theorem \ref{thm.1}  
is to use the Mittag-Leffler
decomposition of the function $q \to 1/(e^q-z)$ (see \cite[Sec.V]{Cn})
\begin{equation}
\lbl{eq.MLexp}
\frac{z}{e^q-z}=  -\frac{1}{2} + \frac{1}{q-\log z}
+ \sum_{k=1}^\infty \frac{1}{q-\log z+2 \pi i k} 
+ \frac{1}{q-\log z-2 \pi i k},
\end{equation}
interchange summation and integration in \eqref{eq.Ia} and use the fact
that
$$
\int_{\ga} q^{\a-1} \frac{dq}{q-\log z-2 \pi i k}
= 2 \pi i (\log z+2 \pi i k)^{\a-1}
$$
\end{remark}

\begin{remark}
\lbl{rem.1}
In \'Ecalle's language, \eqref{eq.link} is a special case of
\begin{equation}
\lbl{eq.major}
\int_0^\infty f(q) dq=
\frac{1}{1-e^{2 \pi i \a}} 
\int_{\ga} \check f(q) dq
\end{equation}
together with the fact that if $f(q)=q^{\a-1}$, then $\check f(q)=
(1-e^{2 \pi i \a})^{-1} f(q)$. For a self-contained introduction to 
majors/minors, see also \cite{Ma}.
\end{remark}

\begin{remark}
\lbl{rem.ML}
The Mittag-Leffler type decomposition \eqref{eq.ML} implies
that $\Li_{\a}(z)$ is multivalued on $\BC\setminus\{0,1\}$, 
for all $\a$ such that $\Re(\a) <0$, $\a\notin\BZ$. It also
implies Equations \eqref{eq.m0} and \eqref{eq.m1}. Indeed, for $z$ near $1$
and $k \neq 0$, $M_{\a}(z)$ is analytic. Thus, \eqref{eq.ML} implies that
for $z$ near $1$ we have
$$
\Li_{\a}(z)=M_{\a}(z) + h(z)
$$
where $h$ is analytic for $z$ near $1$. This proves the second part
of Equation \eqref{eq.m0}. If $z$ is near $0$, then
\begin{eqnarray*}
\Li_{\a}(z \, e^{2 \pi i})-\Li_{\a}(z) &=& M_{\a}[1](z)+
\sum_{k=1}^\infty  M_{\a}[k+1](z) +M_{\a}[-k+1](z) \\ & & -M_{\a}[0](z) -
\sum_{k=1}^\infty  M_{\a}[k](z) +M_{\a}[-k](z)=0.
\end{eqnarray*}
This implies the first part of Equation \eqref{eq.m0}. Equation \eqref{eq.m1}
follows easily from the definition of $M_{\a}(z)$.
\end{remark}

\subsection{Proof of Corollaries \ref{cor.zetas} and \ref{cor.largez}}
\lbl{sub.pfcor}

Corollary \ref{cor.zetas} follows by expanding the sum in \eqref{eq.ML} as a 
convergent power series in $\log z$, and using the functional equation 
for the Riemann zeta function:
$$    
\zeta(s) = 2^s\pi^{s-1}\sin\left(\frac{\pi s}{2}\right)\Gamma(1-s)\zeta(1-s). 
$$
Compare with \cite[Prop.1]{CLZ}.

To prove Corollary \ref{cor.largez}, let us fix $\a$ with $\Re(\a)>0$
and consider the right hand side of Equation \eqref{eq.Llaplace}, which
makes sense for $z\in\BC\setminus[1,\infty)$. The idea is to make some
changes of variables and integration by parts. Let us fix an angle $\th 
\in (0,2 \pi)$ and consider a complex $z=|z| e^{i\th}$ with $|z|$ large.


Making the change of variables $q=\log\tau$,
integrating by parts, and making a change of variables $\tau=z+s$ and
$s=zt$, we obtain that:

\begin{eqnarray*}
z \int_0^\infty q^{\a-1} \frac{z}{e^q-z}
dq &=& z \int_1^\infty \frac{(\log\tau)^{\a-1}}{\tau-z} \frac{d\tau}{\tau}
=
 \frac{z}{\a} \int_1^\infty \frac{1}{\tau-z} d \left((\log\tau)^{\a}\right)
\\
&=&
 \frac{z}{\a} \int_1^\infty \frac{(\log\tau)^{\a}}{(\tau-z)^2} d\tau
=
 \frac{z}{\a} \int_{1-z}^\infty \frac{(\log(z+s))^{\a}}{s^2} d s
\\
&=&
 \frac{1}{\a} \int_{1/z-1}^{\infty e^{-i\th}} \frac{(\log(z+zt))^{\a}}{t^2} dt
=
\frac{1}{\a} \int_{1/z-1}^{\infty e^{-i\th}} \frac{(\log z + \log(t+1))^{\a}}{t^2} dt.
\end{eqnarray*}
Let us separate the domain of integration in two parts: 
$|t| \leq |\log z|$ and $|\log z| \leq |t|$.
The first integral gives:

\begin{eqnarray*}
\frac{1}{\a} \int_{1/z-1}^{\log z} \frac{(\log z + \log(t+1))^{\a}}{t^2} dt
&=&
\frac{(\log z)^{\a}}{\a} 
\int_{1/z-1}^{\log z} \frac{(1 + \log(t+1)/\log z)^{\a}}{t^2} dt.
\end{eqnarray*}
Note that the numerator of the integrand satisfies:

\begin{equation*}
\left(1 + \frac{\log(t+1)}{\log z}\right)^{\a}
=1+O\left(\frac{\log \log z}{\log z}\right)
=1+o(1).
\end{equation*}
So, the first integral gives:

\begin{eqnarray*}
\frac{1}{\a} \int_{1/z-1}^{\log z} \frac{(\log z + \log(t+1))^{\a}}{t^2} dt
&=&\frac{(\log z)^{\a}}{\a} (1+o(1))
\end{eqnarray*}
For the second integral, use 
$
(A+B)^{\a} \leq (2 \max\{A,B\})^{\a} \leq 2^{\a}(A^{\a}+B^{\a})
$
(valid for $A, B \geq 0$ and $\Re(\a)>0$). It follows that we can estimate
the second integral by:

\begin{eqnarray*}
\frac{1}{\a} \int_{\log z}^{\infty e^{-i\th}} 
\left| \frac{(\log z + \log(t+1))^{\a}}{t^2} \right| dt
& \leq &
\frac{2^{\a}}{\a} \int_{\log z}^{\infty e^{-i\th}} \frac{|\log z|^{\a}}{t^2} dt
+
\frac{2^{\a}}{\a} \int_{\log z}^{\infty e^{-i\th}} \frac{|\log(t+1)|^{\a}}{t^2} dt
\\
& = & O\left((\log z)^{\a-1}\right).
\end{eqnarray*}
The result follows.

\begin{remark}
\lbl{rem.transseries}
In fact, we can give a {\em transseries expansion} of 
$\Li_{\a}(z)$ for large $z$ in terms of  $\log z$ and $\log\log z$.
\end{remark}

\ifx\undefined\bysame
        \newcommand{\bysame}{\leavevmode\hbox
to3em{\hrulefill}\,}
\fi

\end{document}